\documentclass[journal,twoside]{IEEEtran}  

\usepackage{graphicx}
\usepackage{stfloats}
\usepackage{tabularx}

\usepackage{epsfig} 
\usepackage{epstopdf}
\usepackage{amssymb}  

\usepackage{cite}
\usepackage{url}
\usepackage{hyperref}
\usepackage[hyphenbreaks]{breakurl}
\usepackage{booktabs}
\usepackage[american,fulldiode]{circuitikz} 

\usepackage{algorithm}
\usepackage[noend]{algpseudocode}
\usepackage{float}
\newfloat{algorithm}{t}{lop}

\usepackage[caption=false,font=footnotesize]{subfig}

\usepackage{tablefootnote}

\usepackage{amsmath} 
\usepackage{amssymb}  

\renewcommand{\leftmark}{{Preprint No: ADP-15-34/T936}}
\renewcommand{\rightmark}{{Preprint No: ADP-15-34/T936}}

\usepackage{tikz}
\usetikzlibrary{arrows}
\usetikzlibrary{automata}

\title{
Toward Topologically Based Upper Bounds on the \\ Number of Power Flow Solutions
}

\author{Daniel K. Molzahn$^{\ast}$, Dhagash Mehta$^{\dagger}$, and Matthew Niemerg$^{\S}$\vspace{-7pt}
\thanks{${\ast}$ Argonne National Laboratory, Energy Systems Division.
        {\tt\small dmolzahn@anl.gov}}%
\thanks{${\dagger}$ University of Notre Dame, Dept. of Applied and Computational Mathematics and Statistics; University of Adelaide, School of Physical Sciences, Dept. of Physics, Centre for the Subatomic Structure of Matter. Support from NSF-ECCS award ID 1509036 and an Australian Research Council DECRA fellowship no. DE140100867.
        {\tt\small dmehta@nd.edu}, \mbox{\scriptsize Preprint No: ADP-15-34/T936}}%
\thanks{${\S}$ Fields Institute for Research in Mathematical Sciences:
        {\tt\small research@matthewniemerg.com}}%
}

\begin{document}

\maketitle


\pagestyle{headings}

\begin{abstract}

The power flow equations, which relate power injections and voltage
phasors, are at the heart of many electric power system computations.
While Newton-based methods typically find the ``high-voltage'' solution to
the power flow equations, which is of primary interest, there are
potentially many ``low-voltage'' solutions that are useful for certain
analyses. This paper addresses the number of solutions to the power flow
equations. There exist upper bounds on the number of power flow solutions;
however, there is only limited work regarding bounds that are functions of
network topology. This paper empirically explores the relationship between
the network topology, as characterized by the maximal cliques, and the
number of power flow solutions. To facilitate this analysis, we use a
numerical polynomial homotopy continuation approach that is guaranteed to
find all complex solutions to the power flow equations. The number of
solutions obtained from this approach upper bounds the number of real
solutions. Testing with many small networks informs the development of
upper bounds that are functions of the network topology. Initial results
include empirically derived expressions for the maximum number of solutions
for certain classes of network topologies.

\end{abstract}

\section{Introduction}

The power flow equations model the steady-state relationship between the voltages and power injections in a power system. Solutions to these nonlinear equations correspond to the equilibria of the underlying differential-algebraic equations describing power system dynamic behavior. Multiple power flow solutions may exist~\cite{klos1975}. Typically, lightly loaded systems have many solutions which disappear in bifurcations with increasing loading. After reaching the system's maximum loadability limit, the power flow equations are infeasible. This paper presents a empirical analysis leading towards development of topologically based upper bounds on the number of power flow solutions.

For typical operating conditions, there is a single ``high-voltage'' solution corresponding to the stable equilibrium of reasonable dynamical models.
The voltage magnitudes associated with this solution are usually close to nominal values, giving rise to the name ``high-voltage'' solution.

There exist mature techniques which often succeed in calculating the high-voltage solution. Typical power flow solvers employ Newton-Raphson and Gauss-Seidel algorithms~\cite{glover_sarma_overbye}. For a sufficiently close initialization, Newton-based techniques yield quadratic convergence to a power flow solution. However, the capabilities of Newton-based techniques are very dependent on the initialization. Reasonable initializations (e.g., the solution to a related problem or a ``flat start'' of $1 \angle 0^\circ$~per~unit voltages) often result in convergence to the high-voltage solution. Obtaining an appropriate initialization is more difficult for parameters outside typical operating ranges, as may occur due to high penetrations of renewable generation and during contingencies. Illustrating this initialization challenge,~\cite{thorp1989} demonstrates that the basins of attraction for Newton-based methods are fractal in nature. 

This challenge motivates the development of alternative solution techniques. Many alternative algorithms are Newton-based with modifications that improve robustness to the choice of initialization. For instance,~\cite{iwamoto1981} describes an ``optimal multiplier'' which modifies the step size to ensure non-divergence of the iterations. 

Many recent approaches do not rely on Newton-based iteration. For instance, the Holomorphic Embedding Load Flow method~\cite{trias2012}, which uses analytic continuation theory from complex analysis, is claimed to be capable of reliably finding a stable solution for any feasible set of power flow equations. A monotone operator approach in~\cite{dj2015} identifies regions with at most one power flow solution. Any solution contained within these regions can be quickly calculated. Convex relaxation techniques can calculate power flow solutions~\cite{lavaei_convexpf} and certify infeasibility~\cite{pfcondition,irep2013,hicss2016}. Progress has also been made using active/reactive power ``decoupling'' approximations~\cite{ilic1992,dorfler2013}.

Many of these approaches focus primarily on calculation of the high-voltage power flow solution. However, other power flow solutions generally exist, often at lower voltages. Locating multiple power flow solutions is important for many types of power system stability assessments~\cite{venikov1975,tamura1983,ribbens-pavella1985,chiang1987}. There may also exist multiple stable power flow solutions, particularly in the presence of power flow reversal conditions on distribution systems~\cite{turitsyn2014}. The literature includes attempts to calculate multiple power flow solutions using a Newton-based algorithm with a range of carefully chosen initializations~\cite{overbye1996,overbye2000}. A semidefinite relaxation is used in \cite{allerton2011} to identify multiple power flow solutions. However, these approaches may not find all solutions.


A numerical continuation approach that claims to find all power flow solutions was presented in~\cite{thorp1993}. Since it scales with the number of \emph{actual} rather than \emph{potential} power flow solutions, this approach is computationally tractable for realistic power systems. However, the robustness proof indicating that this approach finds all solutions to all power flow problems is flawed~\cite{chen2011}, and~\cite{counterexample2013} presents a counterexample. This counterexample also invalidates the robustness of the related approach in~\cite{liu2005} for calculating all power flow solutions that have a certain stability property. Recent work~\cite{lesieutre_wu_allerton2015} presents a related method that improves the robustness of~\cite{thorp1993}.

Although not yet scalable to large systems, there are approaches which are guaranteed to find all power flow solutions. For instance, interval analysis~\cite{mori1999} and Gr\"obner bases~\cite{ning2009} methods are applicable to systems with up to five buses. Other methods guaranteed to find all solutions to systems of polynomial equations are the eigenvalue technique in~\cite{dreesen2009} and the ``moment/sum-of-squares'' relaxations in~\cite{lasserre_book}.

The most computationally tractable methods that reliably find all solutions are based on numerical polynomial homotopy continuation (NPHC). These methods use continuation to trace all complex solutions from a selected polynomial system for which all solutions can be easily characterized to the solutions of the specified target system. In this context, the power flow equations are transformed by splitting real and imaginary parts of the voltage phasors to obtain polynomials in real variables. Thus, only the real solutions to these polynomial equations are physically meaningful.  Nevertheless, ensuring recovery of all real solutions requires a number of continuation traces that depends on an upper bound for the number of complex solutions. Existing techniques are tractable for systems with up to 14 buses~\cite{mehta2014a,mehta2014b} (and the equivalent of $18$ buses for the related Kuramoto model~\cite{mehta2015algebraic}).

The computational complexity of NPHC is unknown and is related to Smale's 17th open problem of the 21st century~\cite{smale1999}. As such, the complexity of NPHC is considered in terms of the maximal number of complex solutions of a polynomial system in a fixed parametrized family of coefficients. Obtaining tighter upper bounds on the number of possible power flow solutions for a specified network thus has the potential to improve the tractability of NPHC. Upper bounds on the number of power flow solutions are also theoretically interesting in their own right.

After introducing the power flow problem, this paper reviews existing bounds and clarifies inaccurate statements in the literature on the maximum number of power flow solutions. To improve upon existing bounds, this paper then uses NPHC to perform an empirical analysis regarding the number of complex power flow solutions. 
By testing a variety of network topologies, the results inform conjectured network-structure dependent upper bounds.

\section{Overview of the Power Flow Equations}
\label{l:power_flow}

The nonlinear power flow equations model the steady-state relationship between voltages and power injections in a power system. Denote the voltage phasors for a $n$-bus system as $V \in \mathbb{C}^n$. We use rectangular coordinates $V_i = V_{di} + \mathbf{j} V_{qi}$, where $V_d,V_q \in \mathbb{R}^n$ and $\mathbf{j}$ is the imaginary unit, to obtain a system of polynomial equations in real variables $V_d$ and $V_q$.

Let $\mathbf{Y} = \mathbf{G} + \mathbf{j} \mathbf{B}$ represent the network admittance matrix. Using ``active/reactive''
representation of complex power injections $P_{i}+\mathbf{j}Q_{i}$, power balance at bus~$i$ yields
\vspace{-5pt}%
\begin{subequations}
{\small
\label{power_flow_equations}
\begin{align}\nonumber \\[-10pt] 
\label{Prect}
P_i = & V_{di} \sum_{k=1}^n \left( \mathbf{G}_{ik} V_{dk} - \mathbf{B}_{ik} V_{qk} \right) + V_{qi} \sum_{k=1}^n \left( \mathbf{B}_{ik}V_{dk} + \mathbf{G}_{ik}V_{qk} \right) \\[0pt] 
\label{Qrect}
Q_i = & V_{di} \sum_{k=1}^n \left( -\mathbf{B}_{ik}V_{dk} - \mathbf{G}_{ik} V_{qk}\right) \!+\! V_{qi} \sum_{k=1}^n \left( \mathbf{G}_{ik} V_{dk} - \mathbf{B}_{ik} V_{qk}\right)
\end{align}
}
The squared voltage magnitude at bus~$i$ is
{\small
\begin{equation} \label{Vmag} 
\left|V_i\right|^2 = V_{di}^2 + V_{qi}^2
\end{equation}
}
\end{subequations}
\vspace{-10pt}

Each bus is classified as PQ, PV, or slack. PQ buses, which typically correspond to loads and are denoted by the set $\mathcal{PQ}$, treat $P_i$ and $Q_i$ as specified quantities and enforce the active power \eqref{Prect} and reactive power \eqref{Qrect} equations. PV buses, which typically correspond to generators and are denoted by the set $\mathcal{PV}$, enforce \eqref{Prect} and \eqref{Vmag} with specified $P_i$ and $\left|V_i\right|^2$. The associated reactive power $Q_i$ may be computed as an ``output quantity'' via \eqref{Qrect}.\footnote{This paper does not consider reactive-power limited generators.} Finally, a single slack bus is selected with specified $V_i$ (typically chosen such that the reference angle is $0^\circ$; i.e., $V_{qi} = 0$). The set $\mathcal{S}$ denotes the slack bus. The active power $P_i$ and reactive power $Q_i$ at the slack bus are determined from \eqref{Prect} and \eqref{Qrect}.

\section{Review of Bounds on the Number of Power Flow Solutions}
\label{l:lit_review}

Zero active and reactive power injections at a bus correspond to either zero current injection or zero voltage magnitude. Systems composed solely of PQ buses with the exception of a single slack bus can thus achieve at least $2^{\left(n-1\right)}$ real solutions when all power injections at PQ buses are zero. This result has led some publications (e.g.,~\cite{thorp1993,guedes2005,feng2015}) to claim that there exist \emph{at most} $2^{\left(n-1\right)}$ real power flow solutions. However, there are systems with more than $2^{\left(n-1\right)}$ real solutions:~\cite{tavora1972,baillieul1984,klos1991} present three-bus systems with six solutions, which is larger than $2^{\left(n-1\right)} = 4$. Thus, while the approaches in~\cite{thorp1993,guedes2005,feng2015} may find \emph{many} power flow solutions, they do not find \emph{all} solutions for all systems.

B\'{e}zout's theorem with $2n-2$ degree-two polynomials bounds the number of power flow solutions at $2^{\left(2n-2\right)}$~\cite{baillieul1982}. The number of complex solutions for sparse systems of polynomial equations with generic coefficients is given by the BKK bound~\cite{sommese2004advances,bates2013numerically}. For this particular power flow formulation for systems with PQ buses, the BKK bound is the same as the total degree bound from B\'ezout's theorem (i.e., $2^{\left(2n-2\right)}$). However, each unique network topology restricts the power flow equations to a parametrized family in the coefficient space, resulting in a deficiency in the number of solutions and the potential for tighter bounds. By exploiting the structure of the power flow equations for systems without PQ buses,~\cite{baillieul1982} proves a tighter bound of 
\begin{equation}
\label{eq:kappa_n}
\kappa_n = \binom{2n-2}{n-1}
\end{equation}
\noindent where $\binom{\;\cdot\;}{\;\cdot\;}$ is the binomial coefficient function. The number of complex solutions for general systems (i.e., a mix of PV and PQ buses) is also bounded by $\kappa_n$~\cite{guo1994}. 

An upper bound of $2^{\left(n-1\right)} + m$, where $m$ is an unspecified network-dependent parameter, is conjectured in~\cite{klos1991}, but no guidance is provided as to the value of $m$.

Under some technical assumptions,~\cite{chiang2015} shows that the maximum number of certain  (``type-1'') unstable power flow solutions associated with a given stable equilibrium is $2\begin{pmatrix}2\left|\mathcal{PV}\right| \\ \left| \mathcal{PV} \right| \end{pmatrix}$, where $\left| \mathcal{PV}\right|$ is the number of PV buses. However, \cite{chiang2015} does not bound the \emph{total} number of solutions.


The aforementioned bounds are not functions of the network topology. The sole topology-dependent bound is given in~\cite{guo1990}: for a system comprised of certain subnetworks which share a single bus, the number of complex solutions is the product of the number of solutions for each subnetwork. In order to significantly improve on the bound $\kappa_n$ for complete networks, we speculate that developing topologically based bounds is possible for a broader class of networks. The remainder of this paper presents an empirical investigation intended to inform such topologically based bounds.

\section{Overview of the Numerical Polynomial Homotopy Algorithm}
\label{l:homotopy}


With the NPHC method \cite{li2003solving,sommese2005numerical,bates2013numerically}, one first creates a start system that is in the 
same parametrized family as the target system. Then, one easily computes all 
the solutions of the start system and deforms these solutions to all the solutions of the target system.  
This deformation of solutions is done by creating a one-parameter embedding between the start system and 
the target system known as a homotopy.  Under suitable conditions, each isolated solution of the target 
system is connected by a homotopy path to at least one solution of the start system.  
By tracking the paths over the complex numbers, one can be guaranteed to find all of the isolated, 
complex solutions of a polynomial system using NPHC.  Once all the complex solutions are found, 
we then consider only the solutions with a small enough imaginary part and deem these to be real.
For all of our computations, we use Bertini v1.5~\cite{bertini}.

\section{Empirical Results and Discussion}
\label{l:results_discussion}

To discover potential topology-dependent bounds, we used NPHC to calculate all complex solutions 
for many small test cases. In this section, we describe the method employed for generating test cases, 
present the results, and discuss some patterns discovered thus far in the data.

A total of 50,000 systems (10,000 for each size from three to seven buses) were randomly constructed. The systems were created by sampling from a uniform distribution for the number of lines ($n+1$ to $\frac{n^2+n}{2}$), with a topology developed from a random spanning tree~\cite{broder1989} augmented with additional lines whose terminal buses were randomly selected. See the appendix for further test case construction details.

We characterize the power system network topology using its \emph{maximal cliques}. A clique is a subgraph where all nodes are linked to one another (i.e., a \emph{complete} subgraph), and a maximal clique is clique that is not a subgraph of another clique.  The set of maximally sized cliques of a given graph can be computed using the Bron-Kerbosch algorithm~\cite{bron-kerbosch}. Despite being NP-Complete~\cite{karp1972}, calculating the maximal cliques of sparse graphs representing power system networks is tractable (e.g., identifying the maximal cliques of the 2383-bus Polish system~\cite{matpower} only required 35 seconds on a laptop). 

\begin{figure}[t]
\centering
\includegraphics[totalheight=0.33\textheight]{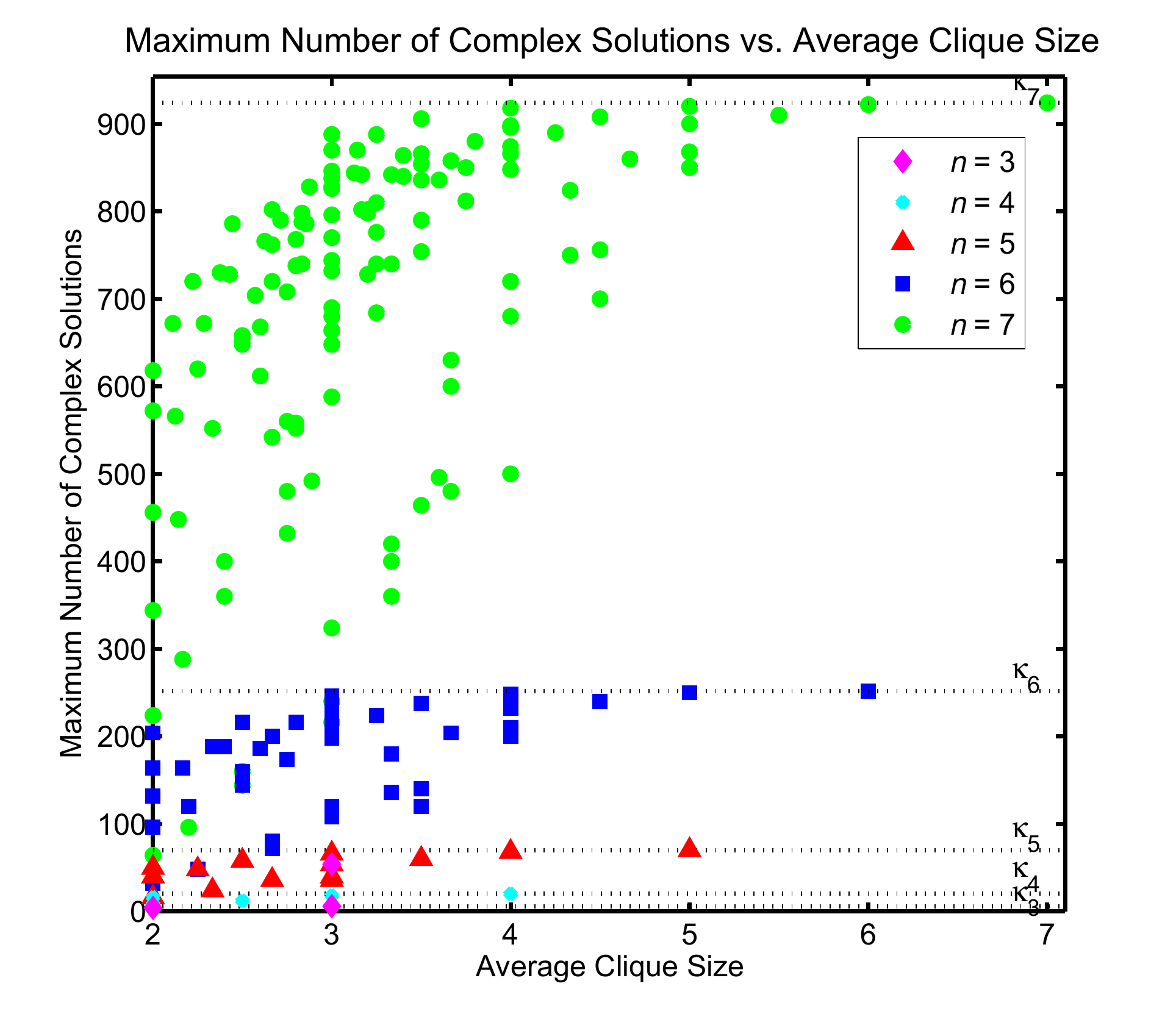}
\vspace{-25pt}
\caption{Maximum number of complex solutions vs. average clique size. Dotted lines correspond to the upper bound $\kappa_n$ for complete networks.}
\vspace{-10pt}
\label{f:numcompsolvsStructSize}
\end{figure}

Fig.~\ref{f:numcompsolvsStructSize} shows the maximum number of complex solutions versus the average clique size (i.e., for a system with maximal cliques $\mathcal{C}_1,\ldots,\mathcal{C}_m$, the average clique size is $\sum_{i=1}^m \left|\mathcal{C}_i\right| / m$, where $\left| \mathcal{C}_i \right|$ is the size of clique~$i$) for each of the 144 unique combinations of clique structures (i.e., unique sets $\{\left|\mathcal{C}_1\right|,\ldots,\left|\mathcal{C}_m\right|\}$) among the 50,000 test cases. Larger average clique size (roughly speaking, a more tightly connected network) is correlated with more complex solutions, although the dependence is not monotonic. 

The dotted lines in Fig.~\ref{f:numcompsolvsStructSize} show the values of $\kappa_n$. The results for the complete networks show that the bounds $\kappa_n$ on the maximum number of complex solutions are exact for systems with up to seven buses. For systems with more than three buses, it is an open question whether there exist power flow equations with a number of \emph{real} solutions equal to $\kappa_n$.\footnote{Two-bus systems can have two real solutions and the three-bus systems in~\cite{tavora1972,baillieul1984,klos1991} have six real solutions, so this question is answered affirmatively for $\kappa_2$ and $\kappa_3$. A complete lossless four-bus network of PV buses with identical line reactances has at most 14 real solutions~\cite{baillieul1982}, which is strictly less than $\kappa_4 = 20$. A variety of perturbations to the electrical parameters for this system still resulted in at most $14$ real solutions.} 

Fig.~\ref{f:numcompsolvsStructSize} illustrates that the number of complex power flow solutions strongly depends on the network structure. The fact that power systems typically have sparse topologies with relatively small maximal cliques suggests that bounds based on the network structure have the potential to be significantly tighter than bounds developed for complete networks.

\begin{figure}[t]
\centering
\vspace{15pt}
\subfloat[Five-bus system with maximal cliques $\mathcal{C}_1 = \{1,2\}$, $\mathcal{C}_2 = \{1,3 \}$, $\mathcal{C}_3 = \{1,4 \}$, and $\mathcal{C}_4 = \{1,5 \}$. Maximum number of solutions is $\kappa^{(1)} = \kappa_2 \kappa_2\kappa_2\kappa_2 = 2\cdot 2\cdot 2 \cdot 2 = 16$.]{\parbox{\columnwidth}{\centering
\begin{tikzpicture}[-,>=stealth',shorten >=1pt,auto,node distance=1cm,
  thick,main node/.style={circle,fill=blue!20,draw,font=\sffamily\bfseries\small,minimum size = 1pt,inner sep = 0.75pt},clique node/.style={rectangle,fill=red!20,draw,font=\sffamily\bfseries\small,minimum size = 1pt,inner sep = 1pt}]

  \node[main node] (1) {1};
  \node[main node] (2) [left of=1] {2};
  \node[main node] (3) [above of=1] {3};
  \node[main node] (4) [right of=1] {4};
  \node[main node] (5) [below of=1] {5};

  \draw (1) edge (2);
  \draw (1) edge (3);
  \draw (1) edge (4);
  \draw (1) edge (5);
  
%
%
\end{tikzpicture}
\label{fig:2x2x2x2_1}}}\\
\subfloat[Four-bus system with maximal cliques $\mathcal{C}_1 = \{1,2,3 \}$ and $\mathcal{C}_2 = \{3,4 \}$. Maximum number of solutions is $\kappa^{(1)} = \kappa_3 \kappa_2 = 2\cdot 6 = 12$.]{\parbox{\columnwidth}{\centering
\begin{tikzpicture}[-,>=stealth',shorten >=1pt,auto,node distance=1.5cm,
  thick,main node/.style={circle,fill=blue!20,draw,font=\sffamily\bfseries\small,minimum size = 1pt,inner sep = 0.75pt}]

  \node[main node] (1) {1};
  \node[main node] (2) [below of=1] {2};
  \node[main node] (3) [right of=2] {3};
  \node[main node] (4) [above of=3] {4};

  \draw (1) edge (2);
  \draw (1) edge (3);
  \draw (2) edge (3);
  \draw (3) edge (4);
\end{tikzpicture}\label{fig:3x2_1}}} \\
%
%
%
%
\subfloat[Seven-bus system with maximal cliques $\mathcal{C}_1 = \{1,2,3 \}$, $\mathcal{C}_2 = \{3,4,5 \}$, and $\mathcal{C}_3 = \{5,6,7 \}$. Maximum number of solutions is $\kappa^{(1)} = \kappa_3 \kappa_3 \kappa_3 = 6\cdot 6\cdot 6 = 216$.]{\parbox{\columnwidth}{\centering
\begin{tikzpicture}[-,>=stealth',shorten >=1pt,auto,node distance=1.5cm,
  thick,main node/.style={circle,fill=blue!20,draw,font=\sffamily\bfseries\small,minimum size = 1pt,inner sep = 0.75pt}]

  \node[draw=none,fill=none] (A) {};
  \node[draw=none,fill=none] (B) [below of=A] {};
  \coordinate (Middle1) at ($(A)!0.5!(B)$);
  
  \node[main node] (3) [right of=Middle1] {3};
  \node[main node] (4) [above right of=3] {4};
  \node[main node] (5) [right of=3] {5};

  \node[main node] (1) [above left of=3] {1};
  \node[main node] (2) [left of=3] {2};

  \node[main node] (7) [right of=4] {7};
  \node[main node] (6) [right of=5] {6};

  \draw (1) edge (2);
  \draw (1) edge (3);
  \draw (2) edge (3);
  \draw (3) edge (4);
  \draw (3) edge (5);
  \draw (4) edge (5);
  \draw (5) edge (7);
  \draw (5) edge (6);
  \draw (6) edge (7);
\end{tikzpicture}
\label{fig:3x3x3_1}}}
\caption{Topologies for example test cases with a block network structure (i.e., the upper bound $\kappa^{(1)}$ in~\eqref{eq:kappa1} is applicable).}
\label{f:kappa1}
\vspace{-10pt}
\end{figure}
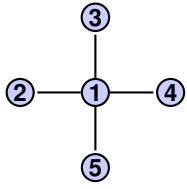
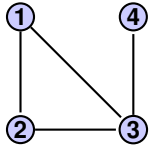
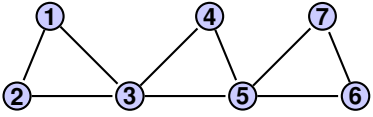

We next describe bounds for two classes of networks. The first class is \emph{block networks}. The blocks of a graph are the maximally sized subgraphs such that the subgraph remains connected after removal of any single node. A network with the property that all blocks are cliques is a block network~\cite{harary1963}. First shown in~\cite{guo1990}, the number of solutions for this class of networks, denoted $\kappa^{\left(1\right)}$, is given by the product of the number of solutions for each clique:
\begin{equation}
\label{eq:kappa1}
\kappa^{\left(1\right)} = \prod_{i=1}^m \kappa_{\left|\mathcal{C}_k \right|}.
\end{equation}
\noindent Fig.~\ref{f:kappa1} shows example block network topologies with corresponding $\kappa^{(1)}$ values. 

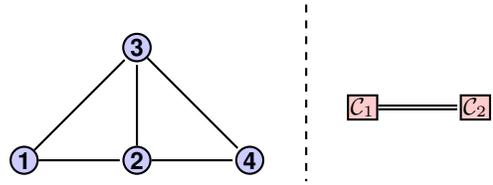
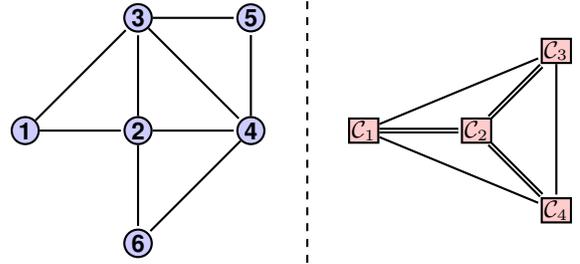
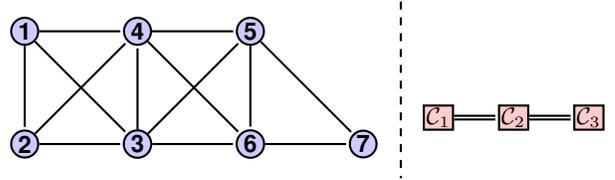
\begin{figure}[t]
\centering
\subfloat[Four-bus system with maximal cliques $\mathcal{C}_1 = \{1,2,3 \}$ and $\mathcal{C}_2 = \{2,3,4 \}$. The (conjectured) maximum number of solutions is $\kappa^{(2)} = \frac{\kappa_3 \kappa_3}{2^{(2-1)}} = \frac{6\cdot 6}{2} = 18$.]{\parbox{\columnwidth}{
\centering
\begin{tikzpicture}[-,>=stealth',shorten >=1pt,auto,node distance=1.5cm,
  thick,main node/.style={circle,fill=blue!20,draw,font=\sffamily\bfseries\small,minimum size = 1pt,inner sep = 0.75pt},clique node/.style={rectangle,fill=red!20,draw,font=\sffamily\bfseries\small,minimum size = 1pt,inner sep = 1pt}]

  \node[main node] (1) {1};
  \node[main node] (2) [right of=1] {2};
  \node[main node] (3) [above of=2] {3};
  \node[main node] (4) [right of=2] {4};
  
  \draw (1) edge (2);
  \draw (1) edge (3);
  \draw (2) edge (3);
  \draw (2) edge (4);
  \draw (3) edge (4);
  
  \node[clique node] (C1) [right of=4, yshift=20pt] {$\mathcal{C}_1$};
  \node[clique node] (C2) [right of=C1] {$\mathcal{C}_2$};

  \draw[double] (C1) -- (C2);
 
  \coordinate (M) at ($(4)!0.5!(C1)$);
  \node[draw=none,fill=none,yshift=10pt] (A) [above of=M] {};
  \node[draw=none,fill=none,yshift=20pt] (B) [below of=M] {};
  \draw[dashed] (A) edge (B);
  
\end{tikzpicture}
\label{fig:3x3_2}}} \\
\subfloat[
Six-bus system with maximal cliques $\mathcal{C}_1 = \{1,2,3 \}$, $\mathcal{C}_2 = \{2,3,4 \}$, $\mathcal{C}_3 = \{3,4,5 \}$, and $\mathcal{C}_4 = \{2,4,6 \}$. The (conjectured) maximum number of solutions is $\kappa^{(2)} = \frac{\kappa_3 \kappa_3 \kappa_3 \kappa_3}{2^{(4-1)}} = \frac{6\cdot 6\cdot 6\cdot 6}{2^3} = 162$.]{\parbox{\columnwidth}{\centering
\begin{tikzpicture}[-,>=stealth',shorten >=1pt,auto,node distance=1.5cm,
  thick,main node/.style={circle,fill=blue!20,draw,font=\sffamily\bfseries\small,minimum size = 1pt,inner sep = 0.75pt},clique node/.style={rectangle,fill=red!20,draw,font=\sffamily\bfseries\small,minimum size = 1pt,inner sep = 1pt}]

  \node[main node] (1) {1};
  \node[main node] (2) [right of=1] {2};
  \node[main node] (3) [above of=2] {3};
  \node[main node] (4) [right of=2] {4};
  \node[main node] (5) [above of=4] {5};  
  
  \draw (1) edge (2);
  \draw (1) edge (3);
  \draw (2) edge (3);
  \draw (2) edge (4);
  \draw (3) edge (4);
  \draw (3) edge (5);
  \draw (4) edge (5);  
  
  \node[main node] (6) [below of=2] {6};
  \draw (4) edge (6);
  \draw (2) edge (6);
  
  \node[clique node] (C1) [right of=4] {$\mathcal{C}_1$};
  \node[clique node] (C2) [right of=C1] {$\mathcal{C}_2$};
  \node[clique node] (C3) [above right of=C2] {$\mathcal{C}_3$};
  \node[clique node] (C4) [below right of=C2] {$\mathcal{C}_4$};
  
  \draw[double] (C1) -- (C2);
  \draw[double] (C2) -- (C3);
  \draw[double] (C2) -- (C4);
  \draw (C1) -- (C3);
  \draw (C1) -- (C4);
  \draw (C3) -- (C4);
 
  \coordinate (M) at ($(4)!0.5!(C1)$);
  \node[draw=none,fill=none,yshift=10pt] (A) [above of=M] {};
  \node[draw=none,fill=none,yshift=-15pt] (B) [below of=M] {};
  \draw[dashed] (A) edge (B);
\end{tikzpicture}
\label{fig:3x3x3_2}}} \\
%
%
\subfloat[Seven-bus system with maximal cliques $\mathcal{C}_1 = \{1,2,3,4 \}$, $\mathcal{C}_2 = \{3,4,5,6 \}$, and $\mathcal{C}_3 = \{5,6,7 \}$. The (conjectured) maximum number of solutions is $\kappa^{(2)} = \frac{\kappa_4 \kappa_4 \kappa_3}{2^{(3-1)}} = \frac{20\cdot 20\cdot 6}{2^2} = 600$.]{\parbox{\columnwidth}{\centering
\begin{tikzpicture}[-,>=stealth',shorten >=1pt,auto,node distance=1.5cm,
  thick,main node/.style={circle,fill=blue!20,draw,font=\sffamily\bfseries\small,minimum size = 1pt,inner sep = 0.75pt},clique node/.style={rectangle,fill=red!20,draw,font=\sffamily\bfseries\small,minimum size = 1pt,inner sep = 1pt,node distance=1cm}]

  \node[main node] (1) {1};
  \node[main node] (2) [below of=1] {2};
  \node[main node] (3) [right of=2] {3};
  \node[main node] (4) [above of=3] {4};
  \node[main node] (5) [right of=4] {5};
  \node[main node] (6) [below of=5] {6};
  \node[main node] (7) [right of=6] {7};
  
  \draw (1) edge (2);
  \draw (1) edge (3);
  \draw (1) edge (4);
  \draw (2) edge (3);
  \draw (2) edge (4);
  \draw (3) edge (4);
  \draw (3) edge (5);
  \draw (3) edge (6);
  \draw (4) edge (5);
  \draw (4) edge (6);
  \draw (5) edge (6);
  \draw (5) edge (7);
  \draw (6) edge (7);
  
  \node[clique node] (C1) [right of=7,yshift=10pt] {$\mathcal{C}_1$};
  \node[clique node] (C2) [right of=C1] {$\mathcal{C}_2$};
  \node[clique node] (C3) [right of=C2] {$\mathcal{C}_3$};
  
  \draw[double] (C1) -- (C2);
  \draw[double] (C2) -- (C3);
 
  \coordinate (M) at ($(7)!0.5!(C1)$);
  \node[draw=none,fill=none,yshift=10pt] (A) [above of=M] {};
  \node[draw=none,fill=none,yshift=20pt] (B) [below of=M] {};
  \draw[dashed] (A) edge (B);
\end{tikzpicture}
\label{fig:4x4x3_2}}}
\caption{Topologies (left) and clique graphs (right) for example test cases where the maximal cliques which share two buses are arranged in a tree (i.e., the conjectured upper bound $\kappa^{(2)}$ in~\eqref{eq:kappa2} is applicable).}
\label{f:kappa2}
\vspace{-10pt}
\end{figure}

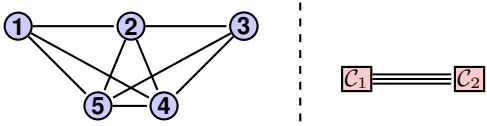
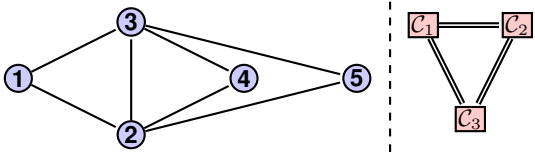
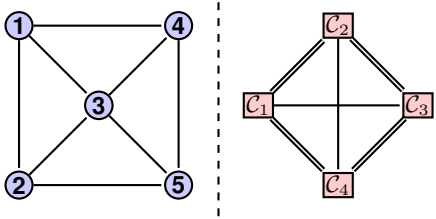
\begin{figure}[!b]
\centering 
\subfloat[Five-bus system with maximal cliques $\mathcal{C}_1 = \{1,2,4,5 \}$ and $\mathcal{C}_2 = \{2,3,4,5 \}$. Maximum of $68$ solutions found among all test cases with isomorphic topologies.]{\parbox{\columnwidth}{\centering

\newlength\triplesep \newlength\triplelinewidth \setlength\triplesep{1pt} \setlength\triplelinewidth{0.75pt} \tikzset{triple/.style={line width=\triplelinewidth,black, preaction={ preaction={draw,line width=2\triplesep+3\triplelinewidth,black}, draw,line width=2\triplesep+\triplelinewidth,white} } }

\begin{tikzpicture}[-,>=stealth',shorten >=1pt,auto,node distance=1.5cm,
  thick,main node/.style={circle,fill=blue!20,draw,font=\sffamily\bfseries\small,minimum size = 1pt,inner sep = 0.75pt},clique node/.style={rectangle,fill=red!20,draw,font=\sffamily\bfseries\small,minimum size = 1pt,inner sep = 1pt}]
  
  \node[main node] (1) {1};
  \node[main node] (2) [right of=1] {2};
  \node[main node] (3) [right of=2] {3};
  \node[main node] (4) [below right of=1] {5};
  \node[main node] (5) [below left of=3] {4};

  \draw (1) edge (2);
  \draw (2) edge (3);
  \draw (2) edge (4);
  \draw (2) edge (5);
  \draw (4) edge (5);
  \draw (1) edge (4);
  \draw (3) edge (5);
  \draw (1) edge (5);
  \draw (3) edge (4);
  
  \node[clique node] (C1) [right of=3,yshift=-20pt] {$\mathcal{C}_1$};
  \node[clique node] (C2) [right of=C1] {$\mathcal{C}_2$};

  \draw[triple] (C1) -- (C2);
  
  \coordinate (M) at ($(3)!0.5!(C1)$);
  \node[draw=none,fill=none,yshift=-20pt] (A) [above of=M] {};
  \node[draw=none,fill=none,yshift=10pt] (B) [below of=M] {};
  \draw[dashed] (A) edge (B);
\end{tikzpicture}
\label{fig:4x4_?}}}\\
\subfloat[Five-bus system with maximal cliques $\mathcal{C}_1 = \{1,2,3 \}$, $\mathcal{C}_2 = \{2,3,4 \}$, and $\mathcal{C}_3 = \{2,3,5 \}$. Maximum of $52$ solutions found among all test cases with isomorphic topologies.]{\parbox{\columnwidth}{\centering
\begin{tikzpicture}[-,>=stealth',shorten >=1pt,auto,node distance=1.5cm,
  thick,main node/.style={circle,fill=blue!20,draw,font=\sffamily\bfseries\small,minimum size = 1pt,inner sep = 0.75pt},clique node/.style={rectangle,fill=red!20,draw,font=\sffamily\bfseries\small,minimum size = 1pt,inner sep = 1pt,node distance=1.25cm}]

  \node[main node] (2) {2};
  \node[main node] (3) [above of=2] {3};
  \coordinate (Middle) at ($(2)!0.5!(3)$);
  \node[main node] (1) [left of=Middle] {1};
  \node[main node] (4) [right of=Middle] {4};
  \node[main node] (5) [right of=4] {5};
  
  \draw (1) edge (2);
  \draw (1) edge (3);
  \draw (2) edge (3);
  \draw (2) edge (4);
  \draw (3) edge (4);
  \draw (3) edge (5);
  \draw (2) edge (5);
  
  \node[clique node] (C1) [above right of=5,yshift=-5pt] {$\mathcal{C}_1$};
  \node[clique node] (C2) [right of=C1] {$\mathcal{C}_2$};
  \coordinate (CM) at ($(C1)!0.5!(C2)$);
  \node[clique node] (C3) [below of=CM] {$\mathcal{C}_3$};
  
  \draw[double] (C1) -- (C2);
  \draw[double] (C1) -- (C3);
  \draw[double] (C2) -- (C3);
  
  \coordinate (M) at ($(5)!0.5!(C1)$);
  \node[draw=none,fill=none,yshift=-20pt] (A) [above of=M] {};
  \node[draw=none,fill=none,yshift=0pt] (B) [below of=M] {};
  \draw[dashed] (A) edge (B);
\end{tikzpicture}
\label{fig:3x3x3_?}}} \\
\subfloat[Five-bus system with maximal cliques $\mathcal{C}_1 = \{1,2,3 \}$, $\mathcal{C}_2 = \{1,3,4 \}$, $\mathcal{C}_3 = \{3,4,5 \}$, and $\mathcal{C}_4 = \{2,3,5 \}$. Maximum of $66$ solutions found among all test cases with isomorphic topologies.]{\parbox{\columnwidth}{\centering
\begin{tikzpicture}[-,>=stealth',shorten >=1pt,auto,node distance=1.5cm,
  thick,main node/.style={circle,fill=blue!20,draw,font=\sffamily\bfseries\small,minimum size = 1pt,inner sep = 0.75pt},clique node/.style={rectangle,fill=red!20,draw,font=\sffamily\bfseries\small,minimum size = 1pt,inner sep = 1pt}]

  \node[draw=none,fill=none] (A) {};
  \node[draw=none,fill=none] (B) [below of=A] {};
  \coordinate (Middle1) at ($(A)!0.5!(B)$);
  
  \node[main node] (3) [right of=Middle] {3};
  \node[main node] (4) [above right of=3] {4};
  \node[main node] (5) [below right of=3] {5};

  \node[main node] (1) [above left of=3] {1};
  \node[main node] (2) [below left of=3] {2};

  \draw (1) edge (2);
  \draw (1) edge (3);
  \draw (2) edge (3);
  \draw (3) edge (4);
  \draw (3) edge (5);
  \draw (4) edge (5);
  \draw (1) edge (4);
  \draw (2) edge (5);
  
  \node[clique node] (C1) [below right of=4] {$\mathcal{C}_1$};
  \node[clique node] (C2) [above right of=C1] {$\mathcal{C}_2$};
  \node[clique node] (C3) [below right of=C2] {$\mathcal{C}_3$};
  \node[clique node] (C4) [below left of=C3] {$\mathcal{C}_4$};
  
  \draw[double] (C1) -- (C2);
  \draw[double] (C1) -- (C4);
  \draw[double] (C2) -- (C3);
  \draw[double] (C3) -- (C4);
  \draw (C1) -- (C3);
  \draw (C2) -- (C4);
    
  \coordinate (M) at ($(4)!0.5!(C1)$);
  \node[draw=none,fill=none,yshift=-15pt] (A) [above of=M] {};
  \node[draw=none,fill=none,yshift=-20pt] (B) [below of=M] {};
  \draw[dashed] (A) edge (B);
\end{tikzpicture}
\label{fig:3x3x3x3_?}}}
\caption{Example test cases without known expressions for upper bounds on the number of solutions.}
\label{f:kappa?}
\end{figure}

The second class of networks is defined based on the corresponding clique graph (i.e., the graph where the maximal cliques are the nodes and the edges represent the buses shared by the maximal cliques). Specifically, the definition of the second class of networks is that the clique graph for those maximal cliques which share exactly two buses (i.e., a single edge) is a tree.\footnote{As illustrated in Figs.~\ref{f:kappa2}\,--\,\ref{f:kappa?}, ``tree'' in this context refers to the lack of loops in connections between maximal cliques which share an edge; i.e., these \emph{maximal cliques}, but not necessarily the \emph{buses}, are connected in a tree.} The NPHC results empirically suggest (but do not prove) that the number of solutions for this second class of networks, denoted $\kappa^{\left(2\right)}$, is at most
\begin{equation}
\label{eq:kappa2}
\kappa^{\left(2\right)} = \frac{1}{2^{m-1}} \prod_{i=1}^m \kappa_{\left|\mathcal{C}_k \right|}.
\end{equation}
\noindent The left hand sides in Fig.~\ref{f:kappa2} show example topologies of systems for which $\kappa^{(2)}$ is applicable. The right hand sides of Fig.~\ref{f:kappa2} show the associated clique graphs, where the nodes represent the maximal cliques and each edge represents a bus shared by the corresponding maximal cliques.

Note that other factors, such as electrical parameters and the number and locations of PV buses, may reduce the actual number of complex solutions for a specific problem below the expression in~\eqref{eq:kappa2}.

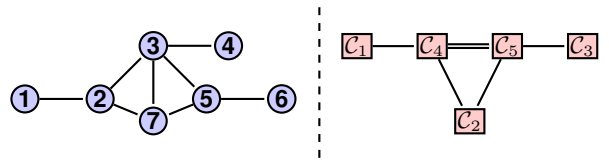
\begin{figure}[!b]
\centering
\begin{tikzpicture}[-,>=stealth',shorten >=1pt,auto,node distance=1cm,
  thick,main node/.style={circle,fill=blue!20,draw,font=\sffamily\bfseries\small,minimum size = 1pt,inner sep = 0.75pt},clique node/.style={rectangle,fill=red!20,draw,font=\sffamily\bfseries\small,minimum size = 1pt,inner sep = 1pt}]

  \node[main node] (1) {1};
  \node[main node] (2) [right of=1] {2};
  \node[main node] (3) [above right of=2] {3};
  \node[main node] (4) [right of=3] {4};
  \node[main node] (5) [below right of=3] {5};
  \node[main node] (6) [right of=5] {6};
  \node[main node] (7) [below of=3] {7};

  \draw (1) edge (2);
  \draw (2) edge (3);
  \draw (2) edge (7);
  \draw (3) edge (4);
  \draw (3) edge (5);
  \draw (3) edge (7);
  \draw (5) edge (6);
  \draw (5) edge (7);
  
  \node[clique node] (C1) [right of=6,yshift=20pt] {$\mathcal{C}_1$};
  \node[clique node] (C4) [right of=C1] {$\mathcal{C}_4$};
  \node[clique node] (C5) [right of=C4] {$\mathcal{C}_5$};
  \coordinate (mid_C4_C5) at ($(C4)!0.5!(C5)$);
  \node[clique node] (C2) [below of=mid_C4_C5] {$\mathcal{C}_2$};
  \node[clique node] (C3) [right of=C5] {$\mathcal{C}_3$};

  \draw (C1) -- (C4);
  \draw[double] (C4) -- (C5);
  \draw (C4) -- (C2);
  \draw (C5) -- (C2);
  \draw (C5) -- (C3);
  
  \coordinate (M) at ($(6)!0.5!(C1)$);
  \node[draw=none,fill=none,yshift=0pt] (A) [above of=M] {};
  \node[draw=none,fill=none,yshift=-10pt] (B) [below of=M] {};
  \draw[dashed] (A) edge (B);
\end{tikzpicture}
\caption{Seven-bus system with maximal cliques $\mathcal{C}_1 = \{1,2\}$, $\mathcal{C}_2 = \{3,4 \}$, $\mathcal{C}_3 = \{5,6\}$, $\mathcal{C}_4 = \{2,3,7\}$, and $\mathcal{C}_5 = \{3,5,7\}$. The (conjectured) maximum number of solutions is $\kappa^{(1)}\kappa^{(2)} = \kappa_2 \kappa_2 \kappa_2 \frac{\kappa_3 \kappa_3}{2^{(2-1)}} = 2\cdot 2\cdot 2 \cdot \frac{6\cdot 6}{2} = 144$. The loop in the clique graph is associated with the maximal cliques sharing a single bus; the restriction of the clique graph to a tree only applies to those cliques sharing \emph{two} buses.}
\label{fig:3x3x2x2x2_1/2}
\end{figure}

To compute the (conjectured) maximum number of solutions for systems where some subnetworks have a block topology and others are arranged such that the maximal cliques have two shared buses (e.g., Fig.~\ref{fig:3x3x2x2x2_1/2}), multiply the values from~\eqref{eq:kappa1} and~\eqref{eq:kappa2} for the appropriate subnetworks. 

The expressions $\kappa^{(1)}$ and $\kappa^{(2)}$ suggest an extension to systems with maximal cliques that share more than two buses: for an integer $\alpha \geq 3$, $\kappa^{\left(\alpha\right)} = \frac{1}{\beta_{\alpha}} \prod_{i=1}^m \kappa_{\left|\mathcal{C}_k \right|}$, where $\beta_{\alpha}$ is some integer function of $\alpha$ (e.g., $\beta_{\alpha} = \alpha^{m-1}$). However, this extension does not hold. See, for instance, the five-bus system in Fig.~\ref{fig:4x4_?}, which has two size-four cliques that share $\alpha = 3$ buses. Since $\kappa_4 = 20$, we would expect that the number of solutions for the five-bus system would have the form $\frac{\kappa_4 \kappa_4}{\beta_3} = \frac{400}{\beta_3}$. The empirical results show a maximum of $68$ solutions for this topology, and $\frac{400}{\beta_3} = 68$ implies that $\beta_3 = 100/17 = 5.8823\ldots$, which is not an integer.

Figs.~\ref{fig:3x3x3_?} and~\ref{fig:3x3x3x3_?} show two other example topologies for which we have not yet been able to determine a pattern in the maximum number of complex solutions. Both of these topologies have a loop in the clique graph among maximal cliques that share two buses. Thus, $\kappa^{(2)}$ is not applicable.

\section{Conclusion}

The power flow equations are intrinsic to many power system analyses. Yet, despite their importance, we lack some basic information such as the number of solutions to these equations for general systems. One practical implication is that numerical polynomial homotopy continuation methods require more computational effort than necessary to find all power flow solutions. After reviewing upper bounds on the number of power flow solutions for complete networks, this paper presented results from a computational experiment using numerical polynomial homotopy continuation to calculate all complex power flow solutions for 50,000 small test cases. These results illustrate the potential for improving existing bounds by exploiting network topology. The results validated the only topology-dependent bound in the existing literature, which is applicable to block networks, and informed the development of conjectured bounds for another class of network topologies.

Future work includes scaling the computations to larger systems, uncovering other patterns to develop further empirical expressions relating the maximum number of power flow solutions to the network topology, and identifying physical explanations for these expressions.




\section*{Appendix}
\label{l:case_details}

For each test case, line resistance $R$, reactance $X$, and shunt susceptance $b$ values for $\Pi$-model equivalent circuits were randomly sampled from Gaussian distributions with mean and standard deviation of $\mu_R = 0.03$, $\sigma_R = 0.03$; $\mu_X = 0.10$, $\sigma_X = 0.03$; and $\mu_b = 0.005$, $\sigma_b = 0.001$ respectively, all in per unit using a $100$~MVA base, with any negative values sampled for line resistances instead set to zero. Lines had an 8\% probability of being transformers with tap ratio $\tau$ and phase-shift $\theta$ sampled from a Gaussian distribution with mean and standard deviation values of $\mu_{\tau} = 1$, $\sigma_{\tau} = 0.02$~per~unit and $\mu_{\theta} = 0^\circ$, $\sigma_{\theta} = 3^\circ$, respectively. A bus was specified to be a generator with 30\% probability, with the first generator selected as the slack bus. If no buses were selected to be generators, a random bus was assigned a generator and chosen to be the slack bus. The generators' voltage setpoints were selected by randomly sampling from a uniform distribution between 0.90 and 1.10~per~unit, and active power injections were sampled from Gaussian distribution with mean and standard deviation values of $\mu_{P_g} = 200$, $\sigma_{P_g} = 30$~MW. (Generator buses had no associated loads.) Buses had a 70\% probability of being specified as loads (PQ buses). Loads have a constant active and reactive power component sampled from a Gaussian distribution with mean and standard deviation $\mu_{P_d} = 50$, $\sigma_{P_d} = 20$~MW and $\mu_{Q_d} = 30$, $\sigma_{Q_d} = 20$~MVAr, respectively, as well as a constant impedance component sampled from a Gaussian distribution with mean and standard deviation $\mu_{P_z} = 0.1$, $\sigma_{P_z} = 0.03$~per~unit and $\mu_{Q_z} = 0.1$, $\sigma_{Q_z} = 0.05$~per~unit for shunt conductance and susceptance, respectively.



\bibliographystyle{IEEEtran}
\bibliography{acc2016}{}

\end{document}